\documentclass[a4paper,12pt]{article}
\usepackage[latin1]{inputenc}
\usepackage{amsmath,amsthm,amsfonts}
\usepackage{amssymb,amstext,amsgen}
\usepackage{amsbsy,amsopn}
\usepackage[pagebackref]{hyperref}
\usepackage{mathrsfs}
\usepackage{amscd}
\usepackage{rotating}
\theoremstyle{plain}
\newtheorem{theorem}{Theorem}[section]

\theoremstyle{definition}

\theoremstyle{remark}

\numberwithin{equation}{section}

\begin{document}

%
%
\newcommand{\R}{\ensuremath{\mathbb R}}
\newcommand{\C}{\ensuremath{\mathbb C}}
\newcommand{\K}{\ensuremath{\mathbb K}}
%
%
\newcommand{\si}{\ensuremath{\sigma}}
\newcommand{\om}{\ensuremath{\omega}}
\newcommand{\Ga}{\ensuremath{\Gamma}}
\newcommand{\Om}{\ensuremath{\Omega}}
%
%
\newcommand{\cc}{\mathcal{C}}
%
%
\newcommand{\Lin}{\ensuremath{{\mathop{\mathrm{L{}}}}}}
\newcommand{\Tan}{\ensuremath{{\mathop{\mathrm{T{}}}}}}
\newcommand{\Cotan}{\ensuremath{{\mathop{\mathrm{T^*}}}}}
%
%
\newcommand{\Dp}{\ensuremath{\mathcal{D}'}}
\newcommand{\DpM}{\ensuremath{\mathcal{D}'(M)}}
\newcommand{\DprsM}{\ensuremath{{\mathcal{D}'}^r_s(M)}}
\newcommand{\DpME}{\ensuremath{{\mathcal{D}'}(M,E)}}
%
%
%
\newcommand{\lgl}{\ensuremath{\langle}}
\newcommand{\rgl}{\ensuremath{\rangle}}
\newcommand{\pr}{{\bf Proof. }}
\newcommand{\ep}{\hspace*{\fill}$\Box$}
\newcommand{\todo}[1]{$\clubsuit$\ {\tt #1}\ $\clubsuit$}
\newcommand{\nn}{\ensuremath{\nonumber}}
\newcommand{\ms}{\medskip\\}
%
%
\newcommand{\TM}{\ensuremath{\mathrm{T}M}}
\newcommand{\TSM}{\ensuremath{\mathrm{T}^*M}}
\newcommand{\TrsM}{\ensuremath{\mathrm{T}^r_s M}}
\newcommand{\cTrsM}{\ensuremath{\mathcal{T}^r_s(M)}}
\newcommand{\cTsrM}{\ensuremath{\mathcal{T}^s_r(M)}}
\newcommand{\CinfM}{\ensuremath{{\mathcal{C}^\infty(M)}}}
\newcommand{\OmncM}{\ensuremath{\Omega^n_c(M)}}
\newcommand{\XM}{\ensuremath{\mathfrak{X}(M)}}
\newcommand{\GaME}{{\ensuremath\Ga(M,E)}}
\newcommand{\GaMF}{{\ensuremath\Ga(M,F)}}
\newcommand{\GaMEs}{{\ensuremath\Ga(M,E^*)}}
\newcommand{\GacME}{{\ensuremath\Ga_c(M,E)}}
\newcommand{\GacMF}{{\ensuremath\Ga_c(M,F)}}
\newcommand{\GacMEs}{{\ensuremath\Ga_c(M,E^*)}}
\newcommand{\GaE}{{\ensuremath\Ga(E)}}
\newcommand{\GaF}{{\ensuremath\Ga(F)}}
\newcommand{\GaEs}{{\ensuremath\Ga(E^*)}}
\newcommand{\GacE}{{\ensuremath\Ga_c(E)}}
\newcommand{\GacF}{{\ensuremath\Ga_c(F)}}
\newcommand{\GacEs}{{\ensuremath\Ga_c(E^*)}}
%
%
%
\newcommand{\DM}{\mathcal{D}(M)}
\newcommand{\Vol}{{\mathop{\mathrm{Vol}}}}
\newcommand{\OmnM}{\Om^n(M)}
\newcommand{\lin}{\mathrm{L}}
\newcommand{\bil}{\mathrm{B}}
\newcommand{\linc}{\mathfrak{L}}
\newcommand{\bilc}{\mathfrak{B}}
\newcommand{\bilsc}{\mathfrak{B}^s}

\title{A Note on Distribution Spaces on Manifolds\footnote{This work was
partially supported by projects P20525 and Y237 of the 
Austrian Science Fund.}}
\author{Michael Grosser\footnote{Faculty of Mathematics,
University of Vienna, Nordbergstra\ss e 15, A-1090, Austria.
e-mail: michael.grosser@univie.ac.at}}
\date{}
\maketitle
\begin{abstract}
26 different concrete representations of the space of vector valued
distributions on a smooth manifold of dimension $n$ are presented
systematically, most of them
new. In the particular case of re\-pre\-sen\-ta\-tions as module homomorphisms
acting on
sections of the dual bundle resp.\ on $n$-forms, the continuity of these
homomorphisms is 
already a consequence of their algebraic properties.
\\[2mm] {\it AMS Mathematics  Subject Classification \kern-1pt(2000)}:
46T30, 46F05, 46H40
\\[1mm] {\it Key words and phrases:} Distributions, differentiable manifolds,
vector bundles, module homomorphisms, automatic continuity
\end{abstract}

\section{Introduction}\label{intro}
For a smooth paracompact orientable Hausdorff manifold $M$ of dimension $n$, the
space $\DpM$ of distributions on $M$ is defined as the topological dual of the
(LF)-space $\OmncM$ of compactly supported $n$-forms on $M$ (for detailed
definitions of these and the following spaces, see Section 3 of \cite{GKOS}). At
many occasions, however, also distributions taking ``values'' in a vector
bundle $E\stackrel{\pi}{\to}M$ (for short: $E$) are required---being aware, of
course, that the notion of the value attained at a point only makes sense for
regular distributions given by, say, continuous functions on $M$. By
\cite{GKOS},
3.1.4 (cf.\ also \cite{CP}, Ch.\ 1, Def.\
7.5; \cite{GS}, Ch.\ VI \S1, p.\ 303; \cite{P}, Section 1; \cite{S},
(1.4.7)) the space of $E$-valued distributions on $M$ is defined as
\begin{equation}\label{DpME}
\DpME:=(\Ga_c(M,E^*\otimes\textstyle{\bigwedge}^n\TSM))'
\end{equation}
where $E^*$ denotes the dual bundle of $E$ and
$\Ga_c(M,E^*\otimes\textstyle{\bigwedge}^n\TSM)$ the space of compactly
supported smooth sections of $E^*\otimes\textstyle{\bigwedge}^n\TSM$
equipped with the usual inductive limit topology. Equation (\ref{DpME}) then
defines $\DpME$ as the topological dual of this latter space of sections. A
brief motivation
for this definition will be given after having presented formula
(\ref{DpMEEE}). Note
that we have replaced $\mathrm{Vol}(M)$ occurring in \cite{GKOS} by
$\bigwedge^n\TSM$ since we are assuming $M$ to be orientable.

In order to make working with $\DpME$ more comfortable various equivalent
representations are used in practice. In a first step, one can rewrite the
section space
$\Ga_c(M,E^*\otimes\textstyle{\bigwedge}^n\TSM)$ as
$\GaMEs\otimes_\CinfM\Ga_c(\bigwedge^n\TSM)=
\GaMEs\otimes_\CinfM\OmncM$, yielding
\begin{equation}\label{DpMEE}
\DpME\cong(\GaMEs\otimes_\CinfM\OmncM)'.
\end{equation}
By (\ref{DpMEE}), we can confine ourselves to
plugging simple tensors $v\otimes\om$ (with $v\in\GaMEs$ and $\om\in\OmncM$) as
arguments into $u\in\DpME$. Moreover, one can show (cf.\
\cite{GKOS}, 3.1.12) that
$(\GaMEs\otimes_\CinfM\OmncM)'$, as a $\CinfM$-module, is isomorphic to the
space
$\Lin_\CinfM(\GaMEs,\DpM)$ of $\CinfM$-linear maps from $\GaMEs$ into $\DpM$,
the correspondence between elements $u\in(\GaMEs\otimes_\CinfM\OmncM)'$ and
$\tilde u\in\Lin_\CinfM(\GaMEs,\DpM)$ given by
$$
\lgl u,v\otimes\om\rgl=\lgl\tilde u(v),\om\rgl
$$
(we always will use angular
brackets to denote the action of any kind of distribution on its respective
arguments). This results in
\begin{equation}\label{DpMEEE}
\DpME\cong\Lin_\CinfM(\GaMEs,\DpM)
\end{equation}
which is also often used in practice.

In order to briefly motivate (\ref{DpMEE}) (and, thereby, Definition
(\ref{DpME})) we start with recalling the definition of (smooth) regular
distributions on $M$. Define $\rho:\CinfM\to\DpM$ by
$\lgl\rho(f),\om\rgl:=\int_M f\om\in\R$ (for $f\in\CinfM$, $\om\in\OmncM$).
Considering, for the sake of simplicity, the vector bundle $E:=\TM$ and the
space $\XM$ of its smooth sections (i.e., of vector fields on $M$), we aim at
obtaining
$\rho:\GaME=\XM\to\Dp(M,\TM)$ by something like
$$
\lgl\rho(X),\fbox{?}\,\rgl:=\int_M X\,\fbox{?}\quad\in\R\qquad (X\in\XM).
$$
Now, for $p\in M$ given, $X(p)$ is a tangent vector at $p$. In order to render a
scalar value for $\int_M X\,\fbox{?}$\,, the sought-after object \fbox{?} must
be
capable of assigning a scalar (depending on $p$, of course, to each $X(p)$ and
then allow integration over $M$.
Thus the most natural choice for \fbox{?} consists in picking a pair
$(\eta,\om)$ consisting of a $1$-form $\eta\in\Om^1(M)=\GaMEs$ and some
$\om\in\OmncM$ and setting $
\lgl\rho(X),(\eta,\om)\rgl:=\int_M (X\cdot\eta)\,\om$ ($X\in\XM$). Noting
that $(X\cdot\eta)\,\om$ is $\CinfM$-bilinear in $\eta$ and $\om$, it is obvious
that $\Om^1(M)\otimes_\CinfM\OmncM$ (rather than $\Om^1(M)\times\OmncM$) is the
appropriate choice for the predual of $\Dp(M,\TM)$. Similarly, for a general
vector bundle $E$ over $M$, we arrive at $\GaMEs\otimes_\CinfM\OmncM$, as
specified at the right hand side of (\ref{DpMEE}). In the trivial case
$E=M\times\R$, this choice also reproduces $\DpM$ as its dual.

The purpose of the present paper is to trace back the algebraic roots of the
relation between (\ref{DpME}), (\ref{DpMEE}) and (\ref{DpMEEE}) and to fully
exploit these structures in order to obtain as many
similar equivalent representations of $\DpME$ as possible. In fact, one also
meets serious topological questions in the course of these investigations:
Observe that $u\in(\GaMEs\otimes_\CinfM\OmncM)'$ is continuous on the whole
space of sections generated by those of the form $v\otimes\om$, whereas for the
corresponding $\tilde u\in\Lin_\CinfM(\GaMEs,\DpM)$, there is no continuity
requirement with respect to $v\in\GaMEs$ whatsoever. Finally, let us mention
that the main impetus towards the present studies was the necessity of
deciding whether $\DprsM:=\Dp(M,\TrsM)$, the space of distributions on $M$
taking values in the bundle of $(r,s)$-tensors, could also be represented as
$\Lin_\CinfM(\OmncM,\cTsrM')$, for $\cTsrM$ denoting the space of smooth
sections of the bundle of $(s,r)$-tensors over $M$. The answer to the latter
question served to solve a problem arising in the intrinsic construction of
generalized tensor fields on smooth manifolds (containing the distributions and
allowing for tensor multiplication and Lie derivatives), see the forthcoming
paper \cite{GKSV}.

Section \ref{mod} will lay the algebraic foundation for the 14 equivalent
representations of the distributions space $\DpME$ to be presented in Section
\ref{topo}. In the latter, the question of linear topologies on
the respective section spaces will be dealt with, allowing to pass from
spaces of linear maps resp.\ linear functionals to spaces of {\it continous}
linear maps
resp.\ {\it continuous} linear functionals, as required by the very definition
of
distribution spaces.
Section \ref{auto}, finally, presents 12 more representations which arise from
automatic continuity. This paper being mainly a report on the results
obtained, we refer to the forthcoming article
\cite{G} for details and complete proofs.

Since throughout the paper, the base manifold of bundles will always be denoted
by $M$, we will, from now on, write section spaces as $\GaE$ etc.\ rather than
$\GaME$ etc.

\section{Module theoretic methods}\label{mod}

The algebraic key to obtaining equivalent representations of $\DpME$ is the
relation
\begin{equation}\label{dualiso}
(V\otimes_A W)^* \cong \bil_A(V,W)\cong
\lin^A(V,W^*)\cong\lin_A(W,V^*).
\end{equation}
Here, $A$ is an algebra over some field $\K$ and $V$, $W$
are linear spaces over $\K$ such that $V$ is a right and $W$ a
left $A$-module. The module tensor product $V\otimes_A W$ is then defined as the
linear space
$V\otimes W/K$ with $K=\mathrm{span}\{(va)\otimes w-v\otimes(aw)\mid a\in A,\
v\in V,\ w\in W\}$. Upper stars next to linear spaces always denoting vector
space duals, $V^*$ and $W^*$ become left resp.\ right [sic!] $A$-modules
by the adjoint action of elements $a\in A$. $\lin^A$ stands for the respective
space of (linear) homorphisms of right $A$-modules, similarly $\lin_A$ for
``left''.
$\bil_A(V,W)$, finally, denotes the space of $A$-balanced
(i.e., satisfying $S(va,w)=S(v,aw)$ for all $a\in A$, $v\in V$, $w\in W$)
bilinear maps from $V\times W$ into $\K$.

The proof of (\ref{dualiso}), based on standard techniques, is straightforward;
compare \cite{Bourb}, Chapter II, \S4.1 for a similar result (based on abelian
groups rather than on linear spaces).
In what follows, we will only consider $\K=\R$ or $\K=\C$.

Setting $A:=\CinfM$, $W:=\GaEs$ and $V:=\OmncM$, (\ref{dualiso})
yields
\begin{eqnarray}
\nn
(\GaEs\otimes_\CinfM \OmncM)^* &\cong& \bil_\CinfM(\GaEs,\OmncM)\\
\label{dualisosect}
&\cong&\lin_\CinfM(\GaEs,\OmncM^*)\\
\nn
&\cong&\lin_\CinfM(\OmncM,\GaEs^*).
\end{eqnarray}
$A=\CinfM$ being commutative in the case at hand, there is no need to
distinguish between left and right modules; therefore we only write
$\lin_\CinfM$.

It is obvious that (\ref{dualisosect}) falls short of giving relevant
information on distribution spaces, due to only the full algebraic duals
occurring. Linear topologies on the section spaces resp.\ on $A$, $V$, $W$
in the general case definitely have to enter the scene. 

Before tackling this question, however, let us introduce another basic
algebraic tool. From now on we generally assume, in addition to the conditions
on $A$, $V$, $W$ stated above, $A$ to be a commutative algebra having an
ideal $B$ with the property
\begin{equation*}
\text{for any given\ }b_1,\ b_2\in B\text{ there exists\ }c\in B\text{
with\ }
cb_1=b_1\text{ and\ }cb_2=b.
\end{equation*}
Denoting the $A$-submodules $BV$ and $BW$ of $V$ resp.\ $W$ by $V_0$ resp.\
$W_0$ one can show that the obvious maps provide canonical $A$-module
homomorphisms
\begin{equation}\label{VWzero}
V_0\otimes_A W \cong V_0\otimes_A W_0\cong V\otimes_A W_0.
\end{equation}
Detailed proofs are given in \cite{G}. Combining relations (\ref{dualiso})
(applied to the pairs $V_0,W$, $V_0,W_0$, $V,W_0$, respectively) and
(\ref{VWzero}), we obtain twelve isomorphic representations of $(V\otimes_A
W_0)^*$.
The latter turns into
$(\GaEs\otimes_\CinfM \OmncM)^*$ by setting $A:=\CinfM$, $B:=\DM$,
$V:=\GaEs$ and $W:=\OmnM$. A corresponding isomorphism result is, of course,
true for general section spaces $V=\GaE$ and $W=\GaF$ where $E$ and $F$ are
arbitrary vector bundles over $M$. We may add $\Ga_c(E\otimes F)^*$ (being
isomorphic to $(\GaE\otimes_\CinfM\GacF)^*$) as the 13th member to this
family.

Finally, noting that $\OmnM$ possesses a non-vanishing section
$\theta$ ($M$ being assumed to be
orientable) we may consider the special case where $W$ is
``faithfully $A$-generated'' in the sense that it is linearly isomorphic to
$A$ via (the inverse of) the map $a\mapsto a\theta$ for some suitable $\theta\in
W$. In \cite{G} it is shown that $V_0\otimes_A W$ is isomorphic to $V_0$ in
this case, allowing to add $\GacEs^*$ (resp.\ $\GacE^*$ in the more general
setting of vector bundles $E,F$) as 14th isomorphic representation to the above
list. Due to the purely algebraic character of the results obtained so far, we
refrain from giving a complete scheme of these spaces, deferring this to the
point where linear topologies have been introduced and successfully dealt with.

\section{Linear topologies}\label{topo}
We endow section spaces $\GaE$ and $\GacE$ with their usual linear Fr\'echet
resp.\ (LF)-topologies (cf., e.g., Section 2 of \cite{GKSV}).
$\GaE\otimes_\CinfM\GaF$ inherits the Fr\'echet topology of $\Ga(E\otimes
F)$ via linear isomorphism and similarly for $\GaE\otimes_\CinfM\GacF$ and
$\Ga_c(E\otimes F)$ with respect to the (LF)-topology. It is easy to see
that the natural map $\tau:\GaE\times
\GaF\to\Ga(E\otimes F)$ sending $(u,v)$ to $u\otimes v = [p\mapsto
u(p)\otimes v(p)]$ is continuous for this choice of topologies. The
corresponding map for compactly supported sections is separately
continuous, at least. Topological duals $X'$ always are endowed with their weak
topologies $\si(X',X)$, and similarly for algebraic duals $X^*$ and
$\si(X^*,X)$.

Now a collection [TOP] of fairly general
conditions 
 on linear
topologies given
on $V$, $V_0$, $W$, $W_0$ and on 
$V_0\otimes_A W$, $V_0\otimes_A W_0$, $V_0\otimes_A W_0$
can be formulated (with every $a\in A$ acting continuously on each of the first
four of them, in particular) allowing, essentially, to replace the notions
``linear'' resp.\
``bilinear'' in the isomorphsm relations of Section \ref{mod} by their
counterparts ``linear contionuous'' resp.\ ``bilinear separately
continuous''. The key to applying this result to section spaces is provided
by the fact that the natural topologies on spaces of sections satisfy conditions
[TOP]. For a detailed account of [TOP], we refer to \cite{G} once more. However,
since we are focusing on the special case of section spaces in this note there
is no need of stating [TOP] explicitly. Rather, we immeditately proceed to the
main theorem given below for
the setting of section spaces: Assuming $F$ to be a line bundle
possessing a non-vanishing section (hence, $\GaF$ to be faithfully
$\CinfM$-generated)
we may pass from the algebraic duals in the isomorphism results of Section
\ref{mod} to the respective topological duals, at the same time replacing
spaces $\lin_\CinfM$ and $\bil_\CinfM$ by the respective spaces $\linc_\CinfM$
of {\it continuous} linear $\CinfM$-module homomorphisms resp.\ $\bilsc_\CinfM$
of {\it separately continuous} bilinear $\CinfM$-balanced maps. Below, the
subscript $\CinfM$ is abbreviated by the subscript $\cc$.

\begin{theorem}\label{alliso}
For vector bundles $E$ and $F$ over $M$ such that $F$ is a line bundle
possessing a non-vanishing section, the topological dual of the
sections space $\Ga_c(E\otimes F)$ with respect to the
(LF)-topology has the following 14
isomorphic representations:
\begin{equation}\label{bigmamaplus}
\begin{aligned}
\Ga&_c(E)'&\ =\ && \Ga&_c(E)'&\ =\ &&
\Ga&_c(E)'\\
&\ \begin{sideways}$\cong$\end{sideways}&&&
&\ \begin{sideways}$\cong$\end{sideways}&&&
&\ \begin{sideways}$\cong$\end{sideways}\\
\Ga_c(E&\otimes F)'&\ =\ && \Ga_c(E&\otimes F)'&\ =\ &&
\Ga_c(E&\otimes F)'\\
&\ \begin{sideways}$\cong$\end{sideways}&&&
&\ \begin{sideways}$\cong$\end{sideways}&&&
&\ \begin{sideways}$\cong$\end{sideways}\\
(\Ga_c(E)&\otimes_{\cc} \Ga(F))' &\ \cong\
&& (\Ga_c(E)&\otimes_{\cc}
\Ga_c(F))' &\ \cong\ && (\Ga(E)&\otimes_{\cc}   \Ga_c(F))'\\
&\ \begin{sideways}$\cong$\end{sideways}&&&
&\ \begin{sideways}$\cong$\end{sideways}&&&
&\ \begin{sideways}$\cong$\end{sideways}\\
\bilsc_{\cc}(\Ga_c&(E),\Ga(F)) &\cong\ &&
\bilsc_{\cc}(\Ga_c&(E),\Ga_c(F))
&\cong\ && \bilsc_{\cc}(\Ga&(E),\Ga_c(F))\\
&\ \begin{sideways}$\cong$\end{sideways}&&&
&\ \begin{sideways}$\cong$\end{sideways}&&&
&\ \begin{sideways}$\cong$\end{sideways}\\
\linc_{\cc}(\Ga_c&(E),\Ga(F)')&\cong\ &&
\linc_{\cc}(\Ga_c&(E),\Ga_c(F)')
&\cong\ && \linc_{\cc}(\Ga(&E),\Ga_c(F)')\\
&\ \begin{sideways}$\cong$\end{sideways}&&&
&\ \begin{sideways}$\cong$\end{sideways}&&&
&\ \begin{sideways}$\cong$\end{sideways}\\
\linc_{\cc}(\Ga(&F),\Ga_c(E)')&\cong\ &&
\linc_{\cc}(\Ga_c&(F),\Ga_c(E)')
&\cong\ && \linc_{\cc}(\Ga_c&(F),\Ga(E)')
\end{aligned}
\end{equation}
\end{theorem}

The complete proof is to be found in \cite{G}.

Replacing $E$ by its dual bundle $E^*$ and setting $F:=\bigwedge^n\TSM$ in the
preceding theorem immediately produces 14 isomorphic representations for the
distribution space $\DpME$.
As a first example, we pick the
isomorphism between $\Ga_c(E\otimes F)'$ and $\linc_\CinfM(\Ga(E),\GacF')$
which, for $E$ replaced by $E^*$ and $F:=\bigwedge^n\TSM$ specializes
to
\begin{equation*}
\DpME\cong\linc_\CinfM(\Ga(E^*),\OmncM')
\end{equation*}
which could be viewed as a ``continuous'' version of the well-known relation
(\ref{DpMEEE}), due to the occurrence of $\linc_\CinfM$ rather than of
$\lin_\CinfM$. As a second example, from
$\Ga_c(E\otimes F)'\cong\linc_\CinfM(\GacF,\GaE')$
we obtain
\begin{equation*}
\DpME\cong\linc_\CinfM(\OmncM,\GaEs')
\end{equation*}
which, similarly, can be considered as a ``continuous'' version of the
relation (stated in Theorem \ref{allisoctd})
$\DpME\cong\lin_\CinfM(\OmncM,\GaEs')$.
In fact, it was the special case $E:=\TrsM$ of the latter
which initiated the present study.

\section{Automatic continuity}\label{auto}
Within the two bottom lines of diagram (\ref{bigmamaplus}), there are 12
(implicit) occurrences of the notion of continuity, altogether: Six in the
form of the symbol $\linc_\CinfM$ denoting {\it continuous} module homomorphism
and six in the form of {\it topological} duals occurring as the range
spaces of the respective $\linc_\CinfM$-spaces. It is a remarkable fact
that under the assumptions of Theorem \ref{alliso}, at all twelve places of
occurrence, the continuity requirement can be omitted (one at a time,
to be precise) without changing the respective space. In other words, for each
of the
twelve instances, an automatic continuity result holds.

Observe that it is certainly not legitimate to drop
two continuity requirements simultaneously, say, by passing from
$\linc_\CinfM(\GaEs),\GacF')$ to $\lin_\CinfM(\GaEs),\GacF^*)$: Replacing
$E$ by $E^*$ and setting
$F:=\bigwedge^n\TSM$, the first of these two spaces, by Theorem \ref{alliso},
represents the distribution space $\DpME$, i.e., the topological dual of
$\GaEs\otimes_\CinfM\OmncM$, whereas the second turns into the {\it
algebraic} dual of $\GaEs\otimes_\CinfM\OmncM$, due to (\ref{dualiso}).

To put proofs of automatic continuity in a nutshell, consider some
left $R$-module homomorphism $T$ of a topological ring $R$ possessing a
right unit $e$
into some left topological $R$-module $X$. Then, for $r\in R$, we have
$T(r)=T(re)=rT(e)$, showing $T$ to be continuous, due to the continuity (with
respect to the first slot) of the action of $R$ on $X$. Needless to say that,
in order to handle the situation at
hand, much more refined methods have to be employed.

It is not feasible in this note to present the proofs of the automatic
continuity results obtained in \cite{G}. As to the six cases of changing
$\linc_\CinfM$ to $\lin_\CinfM$, suffice it to say that for
all three nodes in the last line of (\ref{bigmamaplus}) and for the rightmost
node in the last but one line, direct proofs are feasible, exploiting in some
way or other the fact that $F$ possesses a non-vanishing section. To cover
the remaining positions of the two bottom lines of (\ref{bigmamaplus}), a
transfer principle is employed which holds under some mild assumptions even
in the case of general $A$-modules $V$ and $W$: By this principle, automatic
continuity logically propagates from right to the left in the last but one line,
and from left to right in the last line.

By establishing automatic continuity in the six cases $\linc_\CinfM$ vs.\
$\lin_\CinfM$, we
obtain six more repressentations for $\DpME$. Yet there is more to it: Relation
(\ref{dualiso}) strongly indicates that, in fact, $V$ and $W$ in $L_A(V,W^*)$
enjoy a symmetric position also with respect to automatic continuity. Indeed,
a corresponding statement can be shown to hold for all six nodes of the two
bottom
lines of (\ref{bigmamaplus}) yielding once more a set of six
isomorphic representations
for the space of $E$-valued distributions. In the following concluding Theorem
\ref{allisoctd}, we collect the 14 respresentations obtained by Theorem
\ref{alliso} and the 12 more obtained from automatic continuity as outlined
above. Recall that $\cc$ stands for $\CinfM$.

\begin{theorem}\label{allisoctd}
For vector bundles $E$ and $F$ over $M$ such that $F$ is a line bundle
possessing a non-vanishing section, the topological dual of the
sections space $\Ga_c(E\otimes F)$ with respect to the
(LF)-topology has the following 26
isomorphic representations:
\begin{equation}\label{bigmamactd}
\begin{aligned}
\Ga&_c(E)'&\ =\ && \Ga&_c(E)'&\ =\ &&
\Ga&_c(E)'\\
&\ \begin{sideways}$\cong$\end{sideways}&&&
&\ \begin{sideways}$\cong$\end{sideways}&&&
&\ \begin{sideways}$\cong$\end{sideways}\\
\Ga_c(E&\otimes F)'&\ =\ && \Ga_c(E&\otimes F)'&\ =\ &&
\Ga_c(E&\otimes F)'\\
&\ \begin{sideways}$\cong$\end{sideways}&&&
&\ \begin{sideways}$\cong$\end{sideways}&&&
&\ \begin{sideways}$\cong$\end{sideways}\\
(\Ga_c(E)&\otimes_{\cc} \Ga(F))' &\ \cong\
&& (\Ga_c(E)&\otimes_{\cc}
\Ga_c(F))' &\ \cong\ && (\Ga(E)&\otimes_{\cc}   \Ga_c(F))'\\
&\ \begin{sideways}$\cong$\end{sideways}&&&
&\ \begin{sideways}$\cong$\end{sideways}&&&
&\ \begin{sideways}$\cong$\end{sideways}\\
\bilsc_{\cc}(\Ga_c&(E),\Ga(F)) &\cong\ &&
\bilsc_{\cc}(\Ga_c&(E),\Ga_c(F))
&\cong\ && \bilsc_{\cc}(\Ga&(E),\Ga_c(F))\\
&\ \begin{sideways}$\cong$\end{sideways}&&&
&\ \begin{sideways}$\cong$\end{sideways}&&&
&\ \begin{sideways}$\cong$\end{sideways}\\
\linc_{\cc}(\Ga_c&(E),\Ga(F)')&\cong\ &&
\linc_{\cc}(\Ga_c&(E),\Ga_c(F)')
&\cong\ && \linc_{\cc}(\Ga(&E),\Ga_c(F)')\\
&\ \begin{sideways}$\cong$\end{sideways}&&&
&\ \begin{sideways}$\cong$\end{sideways}&&&
&\ \begin{sideways}$\cong$\end{sideways}\\
\linc_{\cc}(\Ga(&F),\Ga_c(E)')&\cong\ &&
\linc_{\cc}(\Ga_c&(F),\Ga_c(E)')
&\cong\ && \linc_{\cc}(\Ga_c&(F),\Ga(E)')\\
&\ \begin{sideways}$\cong$\end{sideways}&&&
&\ \begin{sideways}$\cong$\end{sideways}&&&
&\ \begin{sideways}$\cong$\end{sideways}\\
\lin_{\cc}(\Ga_c&(E),\Ga(F)')&\cong\ &&
\lin_{\cc}(\Ga_c&(E),\Ga_c(F)')
&\cong\ && \lin_{\cc}(\Ga(&E),\Ga_c(F)')\\
&\ \begin{sideways}$\cong$\end{sideways}&&&
&\ \begin{sideways}$\cong$\end{sideways}&&&
&\ \begin{sideways}$\cong$\end{sideways}\\
\lin_{\cc}(\Ga(&F),\Ga_c(E)')&\cong\ &&
\lin_{\cc}(\Ga_c&(F),\Ga_c(E)')
&\cong\ && \lin_{\cc}(\Ga_c&(F),\Ga(E)')\\
&\ \begin{sideways}$\cong$\end{sideways}&&&
&\ \begin{sideways}$\cong$\end{sideways}&&&
&\ \begin{sideways}$\cong$\end{sideways}\\
\linc_{\cc}(\Ga_c&(E),\Ga(F)^*)&\cong\ &&
\linc_{\cc}(\Ga_c&(E),\Ga_c(F)^*)
&\cong\ && \linc_{\cc}(\Ga(&E),\Ga_c(F)^*)\\
&\ \begin{sideways}$\cong$\end{sideways}&&&
&\ \begin{sideways}$\cong$\end{sideways}&&&
&\ \begin{sideways}$\cong$\end{sideways}\\
\linc_{\cc}(\Ga(&F),\Ga_c(E)^*)&\cong\ &&
\linc_{\cc}(\Ga_c&(F),\Ga_c(E)^*)
&\cong\ && \linc_{\cc}(\Ga_c&(F),\Ga(E)^*)
\end{aligned}
\end{equation}
\end{theorem}

Replacing $E$ by $E^*$ and setting $F:=\bigwedge^n\TSM$ once more we
obtain, from row 7, column 3 of diagram (\ref{bigmamactd}), a new independent
proof of
\begin{equation*}
\DpME\cong\Lin_\CinfM(\GaEs,\DpM)
\end{equation*}
which is (\ref{DpMEEE}) presented in Section \ref{intro} (cf.\ 3.1.12 of
\cite{GKOS}). On the other hand, from row 8, column 3 of 
(\ref{bigmamactd}) we obtain, as a (general) affirmative answer to the
question lying at the origin of this work, the relation
\begin{equation*}
\DpME\cong\Lin_\CinfM(\OmncM,\GaEs').
\end{equation*}
By Theorem \ref{allisoctd}, in either of the preceding relations the
respective spaces of linear module homomorphisms can be replaced by the
corresponding space of {\it continuous} linear module homomorphisms.


\end{document}